\documentclass[twoside,bezier,reqno]{amsart}

\usepackage{amsmath,amsthm,amsfonts,amssymb,amscd,euscript}

\newcommand{\filebegin}{\begin{document}}
\newcommand{\fileend}{\end{document}}
\makeatother
\def\thefootnote{}
\newcommand{\lo}{\longrightarrow}
\newcommand{\NMM}{\hspace*{2mm}}

\renewcommand{\baselinestretch}{1.1}

\renewcommand{\baselinestretch}{1.1}
\def\n{\noindent}%

\numberwithin{equation}{section}
\def\mapdown#1{\Big\downarrow\rlap
{$\vcenter{\hbox{$\scriptstyle#1$}}$}}
\newtheorem{theorem}{Theorem}[section]
\newtheorem{lemma}[theorem]{Lemma}
\newtheorem{proposition}[theorem]{Proposition}
\newtheorem{corollary}[theorem]{Corollary}
\theoremstyle{definition}
\newtheorem{definition}[theorem]{Definition}
\newtheorem{example}[theorem]{Example}
\newtheorem{exercise}[theorem]{Exercise}
\newtheorem{conclusion}[theorem]{Conclusion}
\newtheorem{conjecture}[theorem]{Conjecture}
\newtheorem{criterion}[theorem]{Criterion}
\newtheorem{summary}[theorem]{Summary}
\newtheorem{axiom}[theorem]{Axiom}
\newtheorem{problem}[theorem]{Problem}
\theoremstyle{remark}
\newtheorem{remark}[theorem]{Remark}
\numberwithin{equation}{section}

\begin{document}

\setcounter{page}{1} \noindent

\vspace*{2cm}
\begin{center}
{\bf Controlled continuous $g$-frames and their duals in Hilbert spaces}
 \\[0.5cm]
{Mohamed Rossafi$^{1*}$, Fakhr-dine Nhari$^{2}$ and P. Sam Johnson$^{3}$ \footnote{$^*$Corresponding Author: Mohamed Rossafi} \\[2mm]
$^{1}$LaSMA Laboratory Department of Mathematics Faculty of Sciences, Dhar El Mahraz University Sidi Mohamed Ben Abdellah, Fez, Morocco\\
{\tt E-mail: $^{1}$rossafimohamed@gmail.com$^{*}$}\\
$^{2}$Laboratory Analysis, Geometry and Applications Department of Mathematics, Faculty of Sciences, University of Ibn Tofail, Kenitra, Morocco\\
{\tt E-mail: $^{2}$nharidoc@gmail.com}\\
$^{3}$Department of Mathematical and Computational Sciences, National Institute of Technology Karnataka (NITK), Surathkal, Mangaluru 575 025, India.\\
{\tt E-mail: $^{3}$sam@nitk.edu.in}
} \\[2mm]
\end{center}%
\vspace*{0.5cm}
\begin{quotation}
\noindent
{\footnotesize
{\sc Abstract.}
	In this paper, we characterize and study the concept of  controlled continuous $g$-frame which is an  extension of continuous $g$-frame in Hilbert spaces. We introduce the concept of controlled continuous dual $g$-frame and observe some interesting properties of it. Finally, we characterize all controlled continuous dual $g$-frames of a controlled continuous $g$-frame.}
\end{quotation}
\ \\
{\bf Keywords:} continuous $g$-frame; controlled continuous $g$-frame.\\

\n \textbf{2010 Mathematics subject classification: } 41A58; 42C15.



\section{introduction}

Gabor \cite{Gabor} introduced a method using a family of elementary functions for reconstructing functions (signals)  in $1946$. 
The idea of frames originated in the $1952$ paper by Duffin and Schaeffer  \cite{Duf}  to address some deep questions in non-harmonic Fourier series.  After some decades, Daubechies, Grossmann and Meyer \cite{DGM}  announced formally the definition of frame in the abstract Hilbert spaces in $1986$.  After their work, frame theory began to be widely used, particularly in the more specialized context of wavelet frames and Gabor frames.

Continuous frame is a concept of generalization of frames  proposed by G. Kaiser \cite{KAI} and independently by Ali, Antoine and Gazeau \cite{GAZ} to a family indexed by some locally compact space endowed with a Radon measure. 
These frames are called as frames associated with measurable spaces by Gabrado and Han in \cite{GHN}. Askari-Hemmat, Dehghan and Radjabalipour in \cite{AHDR} called them as generalized frames and in Mathematical Physics they are referred to as coherent states \cite{COH}.

A continuous $g$-frame  is an extension of $g$-frames and continuous frames which was firstly introduced by Abdollahpour and Faroughi in \cite{AFC}. In recent years there has been shown considerable interest by functional analysts in the
study of continuous $g$-frames in Hilbert spaces.
For more details, see \cite{CO, RM}	and the references therein.

We begin with a few preliminaries which are needed in the sequel. Let $H, L$ be separable Hilbert spaces and $(\Omega, \mu)$ be a positive measure space. Let $\{H_{w}\}_{w\in\Omega}$ be a family of closed subspaces of $L$.  Throughout we consider the index set as $\Omega$ and we denote simply $\{H_{w}\}_{w}$ for  $\{H_{w}\}_{w\in\Omega}$.  We denote the set of all bounded linear operators from $H$ into $H_{w}$ by $L(H,H_{w})$  and the set of all bounded linear operators on $H$ with bounded inverse by $GL(H)$.  The set of all positive operators in $GL(H)$ is denoted by $GL^{+}(H)$ and $I_H$ represents the identity operator on $H$.  We note that if $P,Q\in GL(H)$, then $P^{\ast}, P^{-1}$ and $PQ$ are also in $GL(H)$.

\noindent We consider the space
\begin{equation*}
	\ell^{2}(\{H_{w}\}_{w})=\Big\{\{f_{w}\}_{w} : f_{w}\in H_{w},w\in\Omega, \int_{\Omega}\|f_{w}\|^{2}\ d\mu(w)<\infty\Big\}
\end{equation*}
with the inner product given by 
\begin{equation*}
	\langle \{f_{w}\}_{w}, \{g_{w}\}_{w}\rangle =\int_{\Omega}\langle f_{w},g_{w}\rangle \ d\mu(w).
\end{equation*}
It is clear that $\ell^{2}(\{H_{w}\}_{w})$ is a Hilbert space.

\begin{lemma}\cite{CO}\label{lemma_3.1}
	Let $T:H\rightarrow H$ be a linear operator. Then the following are equivalent:
	\begin{itemize}
		\item[(1)] There exist constants $0< A\leq B<\infty$, such that $AI_{H}\leq T\leq BI_{H}$.
		\item[(2)] $T$ is positive and there exist constants $0< A\leq B<\infty$ such that 
		\begin{equation*}
			A\|f\|^{2}\leq \|T^{\frac{1}{2}}f\|^{2}\leq B\|f\|^{2}, \quad \forall f\in H.
		\end{equation*}
		\item[(3)] $T\in GL^{+}(H)$.
	\end{itemize}
\end{lemma}

\begin{definition}\cite{PDA}
	Let $(\Omega,\mu)$ be a measure space with a positive measure $\mu$ and $P\in GL(H)$. A $P$-controlled continuous frame is a map $F:\Omega\rightarrow H$ such that there exist constants $0<A\leq B<\infty$ such that 
	\begin{equation*}
		A\|f\|^{2}\leq \int_{\Omega}\langle f,F(w)\rangle \langle PF(w),f\rangle  \ d\mu(w)\leq B\|f\|^{2},\quad \forall f\in H.
	\end{equation*}
\end{definition}
\begin{definition}\cite{AFC} For each $w\in \Omega$, let $\Lambda_{w}\in L(H, H_w)$. 
	We say that $\Lambda=\{\Lambda_{w}\}_{w}$ is a continuous $g$-frame for $H$ with respect to $\{H_{w}\}_{w}$ if 
	\begin{itemize}
		\item[(1)] for each $f\in H$, $\{\Lambda_{w}f\}_{w}$ is strongly measurable.
		\item[(2)] there exist constants $0<A\leq B<\infty$ such that 
		\begin{equation}\label{eq1}
			A\|f\|^{2}\leq \int_{\Omega}\|\Lambda_{w}f\|^{2}\ d\mu(w)\leq B\|f\|^{2},\quad\forall f\in H.
		\end{equation}
	\end{itemize}
	The numbers $A, B$ are called lower and upper frame bounds for the continuous $g$-frame, respectively. If only the right-hand inequality of \eqref{eq1} is satisfied, we call $\{\Lambda_{w}\}_{w}$
	is a continuous $g$-Bessel family for $H$ with respect to $\{H_{w}\}_{w}$ with bound $B$. If 
	$A =B = \lambda $, we call $\{\Lambda_{w}\}_{w}$ the $\lambda$-tight continuous $g$-frame. Moreover, if $\lambda = 1$, $\{\Lambda_{w}\}_{w}$ is called the Parseval continuous $g$-frame. 
	
	For a given continuous $g$-frame $\Lambda=\{\Lambda_{w}\}_{w}$ for $H$ with respect to $\{H_{w}\}_{w}$, there exists a unique positive and invertible operator (called the frame operator) $S_{\Lambda}:H\to H$ such that for each $f,g\in H$ :
	$$\langle S_\Lambda f, g\rangle=	\int_{\Omega}\langle f,\Lambda_{w}^*\Lambda_{w}g\rangle \ d\mu(w)
	$$ and $AI_{H}\leq S_{\Lambda}\leq BI_{H}$.

\end{definition}
\begin{definition}\cite{YM}
	Let $P,Q\in GL^{+}(H)$ and $\Lambda_{w}\in L(H, H_w)$.  We say that $\{\Lambda_{w}\}_{w}$ is a $(P,Q)$-controlled continuous $g$-frame for $H$ with respect to $\{H_{w}\}_{w}$ if $\{\Lambda_{w}\}_{w}$ is a continuous $g$-Bessel family and there exist constants $0<A\leq B<\infty$ such that 
	\begin{equation}\label{eq3}
		A\|f\|^{2}\leq \int_{\Omega}\langle \Lambda_{w}Pf,\Lambda_{w}Qf\rangle \ d\mu(w)\leq B\|f\|^{2},\quad\forall f\in H.
	\end{equation}
	The numbers $A$ and $B$ are called controlled continuous $g$-frame bounds. If the right hand inequality of \eqref{eq3} holds for all $f\in H$, then $\{\Lambda_{w}\}_{w}$ is a called a $(P,Q)$-controlled continuous $g$-Bessel family with bound $B$.
	
	If $Q=I_{H}$, we call $\{\Lambda_{w}\}_{w}$ is a $P$-controlled continuous $g$-frame for $H$ with respect to $\{H_{w}\}_{w}$. If $Q=P$, we call $\{\Lambda_{w}\}_{w}$ is a $(P,P)$-controlled continuous $g$-frame for $H$ with respect to $\{H_{w}\}$.
\end{definition}

For a $(P,Q)$-controlled continuous $g$-Bessel family $\{\Lambda_{w}\}_{w}$ with bound $B$, the operator $T_{P\Lambda Q}:\ell^{2}(\{H_{w}\}_{w})\rightarrow H$ given by
\begin{equation*}
	T_{P\Lambda Q}\{f_{w}\}_{w}=\int_{\Omega}(PQ)^{\frac{1}{2}}\Lambda_{w}f_{w}\ d\mu(w),\quad \forall \{f_{w}\}_{w}\in \ell^{2}(\{H_{w}\}_{w}) 
\end{equation*}
is well-defined and its adjoint is given by 
\begin{equation*}
	T_{P\Lambda Q}^{\ast}:H\rightarrow \ell^{2}(\{H_{w}\}_{w}),\quad T_{P\Lambda Q}^{\ast}f=\{\Lambda_{w}(QP)^{\frac{1}{2}}f\}_{w},\quad \forall f\in H. 
\end{equation*}
$T_{P\Lambda Q}$ is called the synthesis operator and $T_{P\Lambda Q}^{\ast}$ is called the analysis operator of $\{\Lambda_{w}\}_{w}$. For a $(P,Q)$-controlled continuous $g$-frame $\{\Lambda_{w}\}_{w}$ with bounds $A$ and $B$, the operator 
\begin{equation*}
	S_{P\Lambda Q}:H\rightarrow H, \quad 	S_{P\Lambda Q}f=\int_{\Omega}Q\Lambda_{w}^{\ast}\Lambda_{w}Pf\ d\mu(w),\quad\forall f\in H.
\end{equation*}
is called the frame operator of $\{\Lambda_{w}\}_{w}$. It is a positive invertible operator and $S_{P\Lambda Q}=QS_{\Lambda}P$.

\begin{example}
	Let $P$ and $Q$ be any positive definite matrices of order $2$.  For each $w\in [0, 2\pi]$, we define $\Lambda_{w}=\begin{pmatrix}
		\cos w\\\sin w 
	\end{pmatrix}$.  Then $\{P\Lambda_{w}Q\}$ is a $(P,Q)$-controlled continuous $g$-frame for $\mathbb R^2$.
\end{example}

\section{ controlled continuous $g$-frames in Hilbert spaces}
In this section, we present certain conditions under which continuous $g$-frames become $(P,Q)$-controlled continuous $g$-frames. Given any operators $P,Q\in GL^{+}(H)$, if we have a  $(P,Q)$-controlled continuous $g$-frame $\{\Lambda_{w}\}_{w}$ for $H$ with respect to $\{H_{w}\}_{w}$, we shall get frame bounds of the continuous $g$-frame $\{\Lambda_{w}\}_{w}$ in terms of controlled continuous $g$-frame bounds and operator norms of $(PQ)^{\frac{1}{2}}$ and $(PQ)^{-\frac{1}{2}}$.  It is very useful in estimating bounds of some continuous $g$-frame provided it is a  $(P,Q)$-controlled continuous $g$-frame for some operators $P,Q\in GL^{+}(H)$.

\begin{theorem} Let 	$P,Q\in GL^{+}(H)$. Then $\{\Lambda_{w}\}_{w\in \Omega}$ is a  $(P,Q)$-controlled continuous $g$-frame for $H$ with respect to $\{H_{w}\}_{w}$ if and only if $\{\Lambda_{w}\}_{w}$ is a continuous $g$-frame for $H$ with respect to $\{H_{w}\}_{w}$.
\end{theorem}
\begin{proof}
	Assume that $\{\Lambda_{w}\}_{w}$ is a  $(P,Q)$-controlled continuous $g$-frame for $H$ with respect to $\{H_{w}\}_{w}$ with bounds $A$ and $B$.  Then for each $f\in H$, we have 
	\begin{align*}
		A\|f\|^{2}&=A\|(PQ)^{\frac{1}{2}}(PQ)^{-\frac{1}{2}}f\|^{2}\\&\leq A\|(PQ)^{\frac{1}{2}}\|^{2}\|(PQ)^{-\frac{1}{2}}f\|^{2}\\&\leq\|(PQ)^{\frac{1}{2}}\|^{2}\int_{\Omega}\langle \Lambda_{w}P(PQ)^{-\frac{1}{2}}f,\Lambda_{w}Q(PQ)^{-\frac{1}{2}}f\rangle \ d\mu(w) \\&=\|(PQ)^{\frac{1}{2}}\|^{2}\Big\langle \int_{\Omega}Q\Lambda_{w}^{\ast}\Lambda_{w}P(PQ)^{-\frac{1}{2}}f\ d\mu(w),(PQ)^{-\frac{1}{2}}f\Big\rangle\\&=\|(PQ)^{\frac{1}{2}}\|^{2}\Big\langle QS_{\Lambda}P(PQ)^{-\frac{1}{2}}f,(PQ)^{-\frac{1}{2}}f\Big\rangle\\&=\|(PQ)^{\frac{1}{2}}\|^{2}\langle S_{\Lambda}P^{\frac{1}{2}}Q^{-\frac{1}{2}}f,Q^{\frac{1}{2}}P^{-\frac{1}{2}}f\rangle\\&=\|(PQ)^{\frac{1}{2}}\|^{2}\langle Q^{\frac{1}{2}}P^{-\frac{1}{2}}S_{\Lambda}P^{\frac{1}{2}}Q^{-\frac{1}{2}}f,f\rangle\\&=\|(PQ)^{\frac{1}{2}}\|^{2}\langle S_{\Lambda}f,f\rangle,
	\end{align*}
	hence,
	\begin{equation*}
		\frac{A}{\|(PQ)^{\frac{1}{2}}\|^{2}}\|f\|^{2} \leq \int_{\Omega}\|\Lambda_{w}f\|^{2}\ d\mu(w).
	\end{equation*}
	On the other hand for each $f\in H$, we have 
	\begin{align*}
		\int_{\Omega}\|\Lambda_{w}f\|^{2}\ d\mu(w)&=\langle S_{\Lambda}f,f\rangle \\&=\langle (PQ)^{-\frac{1}{2}}(PQ)^{\frac{1}{2}}S_{\Lambda}f,f\rangle \\&=\langle (PQ)^{\frac{1}{2}}S_{\Lambda}f,(PQ)^{-\frac{1}{2}}f\rangle \\&=\langle S_{\Lambda}(PQ)(PQ)^{-\frac{1}{2}}f,(PQ)^{-\frac{1}{2}}f\rangle \\&=\langle QS_{\Lambda}P(PQ)^{-\frac{1}{2}}f,(PQ)^{-\frac{1}{2}}f\rangle \\&=\langle S_{P\Lambda Q}(PQ)^{-\frac{1}{2}}f,(PQ)^{-\frac{1}{2}}f\rangle\\&\leq B\|(PQ)^{-\frac{1}{2}}\|^{2}\|f\|^{2}.
	\end{align*}
	Finally, we conclude that $\{\Lambda_{w}\}_{w}$ is a continuous $g$-frame for $H$ with respect to $\{H_{w}\}_{w}$ with bounds $\frac{A}{\|(PQ)^{\frac{1}{2}}\|^{2}}$ and $B\|(PQ)^{-\frac{1}{2}}\|^{2}$.

	Conversely, suppose that  $\{\Lambda_{w}\}_{w}$ is a continuous $g$-frame for $H$ with respect to $\{H_{w}\}_{w}$ with bounds $A$ and $B$. Then, 
	\begin{equation*}
		A\langle f,f\rangle \leq \langle S_{\Lambda}f,f\rangle \leq B\langle f,f\rangle, \quad\forall f\in H.
	\end{equation*}
	Since $P,Q\in GL^{+}(H)$, by Lemma \ref{lemma_3.1}, there exist constants $\alpha_{1},\alpha_{2},\beta_{1},\beta_{2}>0$ such that 
	\begin{equation*}
		\alpha_{1}I_{H}\leq P\leq \beta_{1}I_{H}, \quad\alpha_{2}I_{H}\leq Q\leq \beta_{2}I_{H}.
	\end{equation*}
	So, 
	\begin{equation*}
		\alpha_{1}\alpha_{2}A\|f\|^{2}\leq \int_{\Omega}\langle \Lambda_{w}Pf,\Lambda_{w}Qf\rangle \ d\mu(w)\leq \beta_{1}\beta_{2}B\|f\|^{2},  \quad\forall f\in H.
	\end{equation*}
	Therefore  $\{\Lambda_{w}\}_{w}$ is a  $(P,Q)$-controlled continuous $g$-frame for $H$ with respect to $\{H_{w}\}_{w}$.
\end{proof}
\begin{theorem}\label{th1}
	Let $P,Q\in GL^{+}(H)$. Then $\{\Lambda_{w}\}_{w}$ is a  $(P,Q)$-controlled continuous $g$-frame for $H$ with respect to $\{H_{w}\}_{w}$ if and only if $\{\Lambda_{w}\}_{w}$ is a  $((PQ)^{\frac{1}{2}},(PQ)^{\frac{1}{2}})$-controlled continuous $g$-frame for $H$ with respect to $\{H_{w}\}_{w}$.
\end{theorem}	
\begin{proof}
	The proof follows from the following facts :  
	
	\noindent For any $f\in H$, we have 
	\begin{align*}
		\int_{\Omega}\langle \Lambda_{w}Pf,\Lambda_{w}Qf\rangle \ d\mu(w)&=\Big\langle \int_{\Omega}Q\Lambda_{w}^{\ast}\Lambda_{w}Pf\ d\mu(w),f\Big\rangle\\&=\langle QS_{\Lambda}Pf,f\rangle\\&=\langle QPS_{\Lambda}f,f\rangle\\&=\langle(PQ)^{\frac{1}{2}}S_{\Lambda}(PQ)^{\frac{1}{2}}f,f\rangle \\&=\Big\langle \int_{\Omega}(PQ)^{\frac{1}{2}}\Lambda_{w}^{\ast}\Lambda_{w}(PQ)^{\frac{1}{2}}f\ d\mu(w),f\Big\rangle\\&=\int_{\Omega}\langle \Lambda_{w}(PQ)^{\frac{1}{2}}f,\Lambda_{w}(PQ)^{\frac{1}{2}}f\rangle \ d\mu(w).  
	\end{align*}
	Suppose that $\{\Lambda_{w}\}_{w}$ is a  $(P,Q)$-controlled continuous $g$-frame for $H$ with respect to $\{H_{w}\}_{w}$.  Hence
	\begin{equation*}
		A\|f\|^{2}\leq \int_{\Omega}\langle \Lambda_{w}(PQ)^{\frac{1}{2}}f,\Lambda_{w}(PQ)^{\frac{1}{2}}f\rangle \ d\mu(w)\leq B\|f\|^{2},\quad \forall f\in H
	\end{equation*}
	where $A$ and $B$ are continuous frame bounds of $\{\Lambda_{w}\}_{w}$. Thus $\{\Lambda_{w}\}_{w}$ is a  $((PQ)^{\frac{1}{2}},(PQ)^{\frac{1}{2}})$-controlled continuous $g$-frame for $H$ with respect to. $\{H_{w}\}_{w}$.  In a similar way, the other implication can be proved.
\end{proof}
\begin{theorem}
	Let $P,Q\in GL^{+}(H)$. Then $\{\Lambda_{w}\}_{w}$ is a  $(P,Q)$-controlled continuous $g$-frame for $H$ with respect to $\{H_{w}\}_{w}$ if and only if $\{\Lambda_{w}\}_{w}$ is a  $PQ$-controlled continuous $g$-frame for $H$ with respect to $\{H_{w}\}_{w}$.	
\end{theorem}	
\begin{proof}
	The proof is similar to that of Theorem \ref{th1}.
\end{proof}

\begin{definition} Let $(\Omega, \mu)$ be a positive measure space. For each $w\in \Omega$, let  $\Omega_{w}$ be a measurable subset of $\Omega$. We say that the collection 	$\{e_{w,v}\}_{v\in\Omega_{w}}$ is a continuous orthonormal basis of $H_{w}$ if 
	\begin{itemize}
		\item[(1)] $f=\int_{\Omega_{w}}\langle f,e_{w,v}\rangle e_{w,v}\ d\mu(w)$, for all $f\in H_{w}$ ; 
		\item[(2)] $\langle e_{w,v},e_{w,v}\rangle=1$, $\forall v\in \Omega_{w}$ ;
		\item[(3)] $\langle e_{w,v_{1}},e_{w,v_{2}}\rangle =0$, $\forall v_{1},v_{2}\in\Omega_{w}$ and $v_{1}\neq v_{2}$.
	\end{itemize}
\end{definition}
\begin{theorem}\label{th2}
	Let $P,Q\in GL^{+}(H)$. Then $\{\Lambda_{w}\}_{w}$ is a  $(P,Q)$-controlled continuous $g$-frame for $H$ with respect to $\{H_{w}\}_{w}$ if and only if $\{u_{w,v}\}_{w\in\Omega,v\in\Omega_{w}}$ is a  $(P,Q)$-controlled continuous $g$-frame for $H$, where $\{u_{w,v}\}_{w\in\Omega,v\in\Omega_{w}}$ is the sequence induced by $\{\Lambda_{w}\}_{w}$ with respect to $\{e_{w,v}\}_{w\in\Omega,v\in\Omega_{w}}$ (i.e., $u_{w,v}=\Lambda_{w}^{\ast}e_{w,v}$).
\end{theorem}
\begin{proof}
	For each $w\in\Omega$,  let $\{e_{w,v}\}_{v\in\Omega_{w}}$ be an continuous orthonormal basis for $H_{w}$.  Then we have
	\begin{equation*}
		\Lambda_{w}Pf=\int_{\Omega_{w}}\langle \Lambda_{w}Pf,e_{w,v}\rangle  e_{w,v}\ d\mu(v)=\int_{\Omega_{w}}\langle f,P\Lambda_{w}^{\ast}e_{w,v}\rangle e_{w,v}\ d\mu(v)
	\end{equation*}
	and 
	\begin{equation*}
		\Lambda_{w}Qf=\int_{\Omega_{w}}\langle \Lambda_{w}Qf,e_{w,v}\rangle e_{w,v}d\mu(v)=\int_{\Omega_{w}}\langle f,Q\Lambda_{w}^{\ast}e_{w,v}\rangle e_{w,v} \ d\mu(v).
	\end{equation*}
	It is easy to check that 
	\begin{align*}
		\langle \Lambda_{w}Pf,\Lambda_{w}Qf\rangle&=\int_{\Omega_{w}}\langle f,Pu_{w,v}\rangle\langle Qu_{w,v},f\rangle \ d\mu(v)\\
		\int_{\Omega}	\langle \Lambda_{w}Pf,\Lambda_{w}Qf\rangle \ d\mu(w)&=\int_{\Omega}\int_{\Omega_{w}}\langle f,Pu_{w,v}\rangle\langle Qu_{w,v},f\rangle \ d\mu(v)\ d\mu(w).
	\end{align*}
	Hence
	\begin{equation*}
		A\|f\|^{2}\leq\int_{\Omega}\langle \Lambda_{w}Pf,\Lambda_{w}Qf\rangle \ d\mu(w)\leq B\|f\|^{2},\quad\forall f\in H.
	\end{equation*}
	Thus
	\begin{equation*}
		A\|f\|^{2}\leq \int_{\Omega}\int_{\Omega_{w}}\langle f,Pu_{w,v}\rangle\langle Qu_{w,v},f\rangle \ d\mu(v)\ d\mu(w)\leq B\|f\|^{2},\quad\forall f\in H.
	\end{equation*}
	The proof is completed.
\end{proof}
\begin{theorem}
	Let $P,Q\in GL^{+}(H)$. Then $\{\Lambda_{w}\}_{w}$ is a  $(P,Q)$-controlled continuous $g$-frame for $H$ with respect to $\{H_{w}\}_{w}$ if and only if $\{Pu_{w,v}\}_{w\in\Omega,v\in\Omega_{w}}$ is a  $QP^{-1}$-controlled continuous $g$-frame for $H$, where $\{u_{w,v}\}_{w\in\Omega,v\in\Omega_{w}}$ is the sequence induced by $\{\Lambda_{w}\}_{w}$ with respect to $\{e_{w,v}\}_{w\in\Omega,v\in\Omega_{w}}$ (i.e., $u_{w,v}=\Lambda_{w}^{\ast}e_{w,v}$).
\end{theorem}
\begin{proof}
	From the proof of Theorem \ref{th2}, we have 
	\begin{equation*}
		\int_{\Omega}	\langle \Lambda_{w}Pf,\Lambda_{w}Qf\rangle \ d\mu(w)=\int_{\Omega}\int_{\Omega_{w}}\langle f,Pu_{w,v}\rangle\langle Qu_{w,v},f\rangle \ d\mu(v)\ d\mu(w).
	\end{equation*}
	Then 
	\begin{equation*}
		A\|f\|^{2}\leq\int_{\Omega}\langle \Lambda_{w}Pf,\Lambda_{w}Qf\rangle \ d\mu(w)\leq B\|f\|^{2},\quad\forall f\in H,
	\end{equation*}	
	is equivalent to 
	\begin{equation*}
		A\|f\|^{2}\leq \int_{\Omega}\int_{\Omega_{w}}\langle f,Pu_{w,v}\rangle \langle QP^{-1}Pu_{w,v},f\rangle  d\mu(v)\ d\mu(w) \leq B\|f\|^{2}, \quad \forall f\in H.
	\end{equation*}	
	\noindent Hence the proof is completed.
\end{proof}
Combining all results presented in the section, we get the following result which gives several interesting characterizations of  $(P,Q)$-controlled continuous $g$-frames for $H$ with respect to $\{H_{w}\}_{w}$.
\begin{theorem}
	Let $P,Q\in GL^{+}(H)$. Then the following are equivalent:
	\begin{enumerate}
		\item $\{\Lambda_{w}\}_{w}$ is a  $(P,Q)$-controlled continuous $g$-frame for $H$ with respect to $\{H_{w}\}_{w}$.
		\item $\{\Lambda_{w}\}_{w}$ is a  $((PQ)^{\frac{1}{2}},(PQ)^{\frac{1}{2}})$-controlled continuous $g$-frame for $H$ with respect to $\{H_{w}\}_{w}$.
		\item $\{\Lambda_{w}\}_{w}$ is a  $QP$-controlled  continuous $g$-frame for $H$ with respect to $\{H_{w}\}_{w}$.
		\item $\{u_{w,v}\}_{w\in\Omega,v\in\Omega_{w}}$ is a  $(P,Q)$-controlled continuous $g$-frame for $H$, where $\{u_{w,v}\}_{w\in\Omega,v\in\Omega_{w}}$ is the sequence induced by $\{\Lambda_{w}\}_{w}$ with respect to $\{e_{w,v}\}_{w\in\Omega,v\in\Omega_{w}}$.
		\item $\{Pu_{w,v}\}_{w\in\Omega,v\in\Omega_{w}}$ is a  $QP^{-1}$-controlled continuous $g$-frame for $H$, where $\{u_{w,v}\}_{w\in\Omega,v\in\Omega_{w}}$ is the sequence induced by $\{\Lambda_{w}\}_{w}$ with respect to $\{e_{w,v}\}_{w\in\Omega,v\in\Omega_{w}}$.
		\item $\{\Lambda_{w}\}_{w}$ is a continuous $g$-frame for $H$ with respect to$\{H_{w}\}_{w}$.
	\end{enumerate}
	
	Moreover, $\{\Lambda_{w}\}_{w}$ is a continuous $g$-frame for $H$ with respect to $\{H_{w}\}_{w}$ with bounds $\frac{A}{\|(PQ)^{\frac{1}{2}}\|^{2}}$ and $B\|(PQ)^{-\frac{1}{2}}\|^{2}$.
\end{theorem}
\section{Controlled continuous dual $g$-frames in Hilbert spaces}
In this section, we define  controlled continuous dual $g$-frames and characterize them by operator theory.
\begin{definition}\label{def1}
	Let $P,Q\in GL^{+}(H)$ and let $\{\Lambda_{w}\}_{w}, \{\Gamma_{w}\}_{w}$ be  $(P,P)$-controlled continuous and  $(Q,Q)$-controlled continuous $g$-Bessel families for $H$ with respect to $\{H_{w}\}_{w}$, respectively. 
	We say that  $\{\Gamma_{w}\}_{w}$ is  a  $(P,Q)$-controlled continuous dual $g$-frame of  $\{\Lambda_{w}\}_{w}$ if \begin{equation*}
		f=\int_{\Omega}P\Lambda_{w}^{\ast}\Gamma_{w}Qf\ d\mu(w), \quad \forall f\in H.
	\end{equation*}
	In particular, if $Q=I_{H}$,  $\{\Gamma_{w}\}_{w}$ is called a  $P$-controlled continuous dual $g$-frame of $\{\Lambda_{w}\}_{w}$.
\end{definition} 
\begin{definition}
	Let  $P,Q\in GL^{+}(H)$ and let $\{\Lambda_{w}\}_{w}, \{\Gamma_{w}\}_{w}$ be  $(P,P)$-controlled continuous and  $(Q,Q)$-controlled continuous $g$-Bessel families for $H$ with respect to $\{H_{w}\}_{w}$ respectively. 
	For the pair $\{\Lambda_{w}\}_{w}$ and $\{\Gamma_{w}\}_{w}$, we define a  $(P,Q)$-controlled continuous dual $g$-frame operator $S_{P\Lambda\Gamma Q}$ by	\begin{equation*}
		S_{P\Lambda\Gamma Q}f=\int_{\Omega}P\Lambda_{w}^{\ast}\Gamma_{w}Qf\ d\mu(w),\quad\forall f\in H.
	\end{equation*} 
	It is easy to check that $S_{P\Lambda\Gamma Q}$ is a well-defined bounded operator, and 
	\begin{equation*}
		S_{P\Lambda\Gamma Q}=T_{P\Lambda P}T_{Q\Gamma Q}^{\ast}=PT_{\Lambda}T_{\Gamma}^{\ast}Q=PS_{\Lambda\Gamma}Q,
	\end{equation*}
	where $S_{\Lambda\Gamma}f=\int_{\Omega}\Lambda_{w}^{\ast}\Gamma_{w}f\ d\mu(w)$. Note that  $\{\Gamma_{w}\}_{w}$ is a  $(P,Q)$-controlled continuous dual $g$-frame of $\{\Lambda_{w}\}_{w}$ if and only if $S_{P\Lambda\Gamma Q}=I_{H}$.
\end{definition}
\begin{theorem}\label{th3}
	Let $P,Q\in GL^{+}(H)$ and let $\{\Lambda_{w}\}_{w}, \{\Gamma_{w}\}_{w}$ be  $(P,P)$-controlled continuous and  $(Q,Q)$-controlled continuous $g$-Bessel families with bounds $B_{\Lambda}$ and $B_{\Gamma}$ respectively. If $S_{P\Lambda\Gamma Q}$ is bounded below, then  $\{\Lambda_{w}\}_{w}$ and  $\{\Gamma_{w}\}_{w}$ are  $(P,P)$-controlled continuous and  $(Q,Q)$-controlled continuous $g$-frames respectively.
\end{theorem}
\begin{proof}
	Suppose that there exists a constant $\lambda>0$ such that 
	\begin{equation*}
		\|S_{P\Lambda\Gamma Q}f\|\geq \lambda \|f\|, \quad \forall f\in H.
	\end{equation*}
	Hence 
	\begin{align*}
		\lambda\|f\|&\leq \|S_{P\Lambda\Gamma Q}f\|\\&=\sup_{\|g\|=1}\Big|\langle \int_{\Omega}P\Lambda_{w}^{\ast}\Gamma_{w}Qf\ d\mu(w),g\rangle\Big |\\&=\sup_{\|g\|=1}\Big|\int_{\Omega}\langle \Gamma_{w}Qf,\Lambda_{w}Pg\rangle \ d\mu(w) \Big|\\&\leq \sup_{\|g\|=1}\bigg(\int_{\Omega}\|\Gamma_{w}Qf\|^{2}\ d\mu(w))\bigg)^{\frac{1}{2}}\bigg(\int_{\Omega}\|\Lambda_{w}Pg\|^{2}\ d\mu(w)\bigg)^{\frac{1}{2}}\\&\leq \sqrt{B_{\Lambda}}\bigg(\int_{\Omega}\|\Gamma_{w}Qf\|^{2}\ d\mu(w)\bigg)^{\frac{1}{2}}.
	\end{align*}
	Thus 
	\begin{equation*}
		\frac{\lambda^{2}}{B_{\Lambda}}\|f\|^{2}\leq \int_{\Omega}\|\Gamma_{w}Qf\|^{2}\ d\mu(w),\quad \forall f\in H.
	\end{equation*}
	On the other hand, Since 
	\begin{equation*}
		S_{P\Lambda\Gamma Q}^{\ast}=(PS_{\Lambda\Gamma}Q)^{\ast}=QS_{\Gamma\Lambda}P=S_{Q\Gamma\Lambda P},
	\end{equation*}
	then $S_{Q\Gamma\Lambda P}$ is also bounded below. Similarly, we can prove that $\{\Lambda_{w}\}_{w}$ is a  $(P,P)$-controlled continuous $g$-frame. Hence the proof is completed.
\end{proof}
\begin{theorem}
	Let $P,Q\in GL^{+}(H)$, $\{\Lambda_{w}\}_{w}$ and let $\{\Gamma_{w}\}_{w}$ be  $(P,P)$-controlled continuous and  $(Q,Q)$-controlled continuous $g$-Bessel families for $H$ with respect to $\{H_{w}\}_{w}$ respectively. Then the following conditions are equivalent:
	\begin{itemize}
		\item[(1)] $f=\int_{\Omega}P\Lambda_{w}^{\ast}\Gamma_{w}Qf\ d\mu(w)$,$\ \forall f\in H$.
		\item[(2)] $f=\int_{\Omega}Q\Gamma_{w}^{\ast}\Lambda_{w}Pf\ d\mu(w)$,$\ \forall f\in H$.
		\item[(3)] $\langle f,g\rangle=\int_{\Omega}\langle \Lambda_{w}Pf,\Gamma_{w}Qg\rangle \ d\mu(w)=\int_{\Omega}\langle \Gamma_{w}Qf,\Lambda_{w}Pg\rangle \ d\mu(w), \  \forall f,g\in H$.
		\item[(4)] $\|f\|^{2}=\int_{\Omega}\langle \Lambda_{w}Pf,\Gamma_{w}Qf\rangle \ d\mu(w)=\int_{\Omega}\langle \Gamma_{w}Qf,\Lambda_{w}Pf\rangle \ d\mu(w),\ \forall f\in H$.
	\end{itemize}
	In case the equivalent conditions are satisfied, $\{\Lambda_{w}\}_{w}$ and $\{\Gamma_{w}\}_{w}$ are  $(P,P)$-controlled continuous and  $(Q,Q)$-controlled continuous $g$-frames respectively.
\end{theorem}
\begin{proof}
	$(1)\Longleftrightarrow (2)$. Let $T_{P\Lambda P}$ be the synthesis operator of the  $(P,P)$-controlled continuous $g$-Bessel family $\{\Lambda_{w}\}_{w}$ and $T_{Q\Gamma Q}$ be  the synthesis operator of the  $(P,P)$-controlled continuous $g$-Bessel family $\{\Gamma_{w}\}_{w}$. Condition $(1)$ means that $T_{P\Lambda P}T_{Q\Gamma Q}^{\ast}=I_{H}$, and it  is equivalent to $T_{Q\Gamma Q}T_{P\Lambda P}^{\ast}$, which is identical to the statement $(2)$. Conversely, $(2)\implies (1)$ similarly.
	
	$(2)\Longleftrightarrow (3)$. It is clear that $(2)\implies (3)$. For any $f,g\in H$, $\langle f,g\rangle=\int_{\Omega}\langle \Lambda_{w}Pf,\Gamma_{w}Qg\rangle \ d\mu(w)$ shows that 
	\begin{equation*}
		\Big\langle f-\int_{\Omega}Q\Gamma_{w}^{\ast}\Lambda_{w}Pf\ d\mu(w),g\Big\rangle=0,\quad\forall g\in H.
	\end{equation*}
	Hence $(3)\Longrightarrow (2)$ is proved.
	
	$(3)\Longleftrightarrow (4)$. $(3)\implies (4)$ is obvious. To prove $(4)\implies (3)$, we apply condition $(4)$ and get
	\begin{align*}
		\|f+g\|^{2}&=\int_{\Omega}\langle \Lambda_{w}P(f+g),\Gamma_{w}Q(f+g)\rangle \ d\mu(w)\\&=\int_{\Omega}\langle \Lambda_{w}Pf,\Gamma_{w}Qf\rangle \ d\mu(w)+\int_{\Omega}\langle \Lambda_{w}Pf,\Gamma_{w}Qg\rangle \ d\mu(w)\\&\quad +\int_{\Omega}\langle \Lambda_{w}Pg,\Gamma_{w}Qf\rangle \ d\mu(w)+\int_{\Omega}\langle \Lambda_{w}Pg,\Gamma_{w}Qg\rangle \ d\mu(w).
	\end{align*}
	Similarly,
	\begin{align*}
		\|f-g\|^{2}&=\int_{\Omega}\langle \Lambda_{w}Pf,\Gamma_{w}Qf\rangle \ d\mu(w)-\int_{\Omega}\langle \Lambda_{w}Pf,\Gamma_{w}Qg\rangle \ d\mu(w)\\&\quad -\int_{\Omega}\langle \Lambda_{w}Pg,\Gamma_{w}Qf\rangle \ d\mu(w)+\int_{\Omega}\langle \Lambda_{w}Pg,\Gamma_{w}Qg\rangle \ d\mu(w).
	\end{align*}
	\begin{align*}
		\|f+ig\|^{2}&=\int_{\Omega}\langle \Lambda_{w}Pf,\Gamma_{w}Qf\rangle \ d\mu(w)-i\int_{\Omega}\langle \Lambda_{w}Pf,\Gamma_{w}Qg\rangle \ d\mu(w)\\&\quad +i\int_{\Omega}\langle \Lambda_{w}Pg,\Gamma_{w}Qf\rangle \ d\mu(w)+\int_{\Omega}\langle \Lambda_{w}Pg,\Gamma_{w}Qg\rangle \ d\mu(w).
	\end{align*}
	\begin{align*}
		\|f-ig\|^{2}&=\int_{\Omega}\langle \Lambda_{w}Pf,\Gamma_{w}Qf\rangle \ d\mu(w)+i\int_{\Omega}\langle \Lambda_{w}Pf,\Gamma_{w}Qg\rangle \ d\mu(w)\\&\quad -i\int_{\Omega}\langle \Lambda_{w}Pg,\Gamma_{w}Qf\rangle \ d\mu(w)+\int_{\Omega}\langle \Lambda_{w}Pg,\Gamma_{w}Qg\rangle \ d\mu(w).
	\end{align*}
	By polarization identity,
	\begin{align*}
		\langle f,g\rangle&=\frac{1}{4}\bigg(\|f+g\|^{2}-\|f-g\|^{2}+i\|f+ig\|^{2}-i\|f-ig\|^{2}\bigg)\\&=\int_{\Omega}\langle \Lambda_{w}Pf,\Gamma_{w}Qg\rangle \ d\mu(w).
	\end{align*}
	In case the equivalent conditions are satisfied, $S_{Q\Gamma\Lambda P}=I_{H}$ implies $\|S_{Q\Gamma\Lambda P}\|=1$, hence by Theorem \ref{th3}, $\{\Lambda_{w}\}_{w}$ and $\{\Gamma_{w}\}_{w}$ are  $(P,P)$-controlled continuous and  $(Q,Q)$-controlled continuous $g$-frames respectively.
\end{proof}
\begin{theorem}
	Let $P,Q\in GL^{+}(H)$. A sequence $\{\Lambda_{w}\}_{w}$ is a  $(P,Q)$-controlled continuous $g$-Bessel family for $H$ with respect to $\{H_{w}\}_{w}$ with bound $B$ if and only if the operator $T_{P\Lambda Q}:\ell^{2}(\{H_{w}\}_{w})\rightarrow H$ given by 
	\begin{equation*}
		T_{P\Lambda Q}(\{f_{w}\}_{w})=\int_{\Omega}(PQ)^{\frac{1}{2}}\Lambda_{w}^{\ast}f_{w}\ d\mu(w)
	\end{equation*}
	is well-defined and bounded with $\|T_{P\Lambda Q}\|\leq \sqrt{B}$.
\end{theorem}
\begin{proof}
	The necessary condition follows from the definition of  $(P,Q)$-controlled continuous $g$-Bessel sequence. We only need to prove the sufficient condition. Suppose that $T_{P\Lambda Q}$ is well-defined and bound operator with $\|T_{P\Lambda Q}\|\leq \sqrt{B}$. For any $f\in H$, we have 
	\begin{align*}
		\int_{\Omega}\langle \Lambda_{w}Pf,\Lambda_{w}Qf\rangle \ d\mu(w)&=\int_{\Omega}\langle Q\Lambda_{w}^{\ast}\Lambda_{w}Pf,f\rangle \ d\mu(w) \\&=\langle QS_{\Lambda}Pf,f\rangle \\&=\langle (QP)^{\frac{1}{2}}S_{\Lambda}(QP)^{\frac{1}{2}}f,f\rangle\\&=\Big\langle \int_{\Omega}(QP)^{\frac{1}{2}}\Lambda_{w}^{\ast}\Lambda_{w}(QP)^{\frac{1}{2}}f \ d\mu(w),f\Big\rangle \\&\leq \|T_{P\Lambda Q}\|\bigg(\int_{\Omega}\|\Lambda_{w}(QP)^{\frac{1}{2}}f\|^{2}\ d\mu(w)\bigg)^{\frac{1}{2}}\|f\|\\&=\|T_{P\Lambda Q}\|\bigg(\int_{\Omega}\langle \Lambda_{w}Pf,\Lambda_{w}Qf\rangle \ d\mu(w)\bigg)^{\frac{1}{2}}\|f\|.
	\end{align*}
	Hence we get 
	\begin{equation*}
		\int_{\Omega}\langle  \Lambda_{w}Pf,\Lambda_{w}Qf\rangle \ d\mu(w)\leq  \|T_{P\Lambda Q}\|^{2}\|f\|^{2}\leq B\|f\|^{2}.
	\end{equation*}
	This shows that $\{\Lambda_{w}\}_{w}$ is a  $(P,Q)$-controlled continuous $g$-Bessel family for $H$ with respect to $\{H_{w}\}_{w}$ with bound $B$.
\end{proof}
\begin{theorem}
	Let $P,Q\in GL^{+}(H)$ and  $\{\Lambda_{w}\}_{w}$ be a $(P,P)$-controlled $g$-frame for $H$ with respect to $\{H_{w}\}_{w}$ with the synthesis operator $T_{P\Lambda P}$. Then a  $(Q,Q)$-controlled continuous $g$-frame $\{\Gamma_{w}\}_{w}$ is a  $(P,Q)$-controlled continuous dual $g$-frame of $\{\Lambda_{w}\}_{w}$ if and only if 
	\begin{equation*}
		Q\Gamma_{w}^{\ast}e_{w,v}=U(e_{w,v}\delta_{w}),\quad w\in\Omega, v\in\Omega_{w}
	\end{equation*} 
	where $U:\ell^{2}(\{H_{w}\}_{w})\rightarrow H$ is a bounded left-inverse of $T_{P\Lambda P}^{\ast}$.
\end{theorem}
\begin{proof}
	If $\{g_{w}\}_{w}\in \ell^{2}(\{H_{w}\}_{w})$, then 
	\begin{equation*}
		\{g_{w}\}_{w}=\int_{\Omega}g_{w}\delta_{w}\ d\mu(w)=\int_{\Omega}\int_{\Omega_{w}}\langle g_{w},e_{w,v}\rangle  e_{w,v}\delta_{w} \ \ d\mu(v)\ d\mu(w).
	\end{equation*}
	Roughly speaking $\{e_{w,v}\delta_{w}\}_{w\in\Omega, v\in\Omega_{w}}$ is a continuous orthonormal basis of $ \ell^{2}(\{H_{w}\}_{w})$.
	
	If there exists $U:\ell^{2}(\{H_{w}\}_{w})\rightarrow H$ is a bounded left-inverse of $T_{P\Lambda P}^{\ast}$ such that 
	\begin{equation*}
		Q\Gamma_{w}^{\ast}e_{w,v}=U(e_{w,v}\delta_{w}),\quad w\in\Omega, v\in\Omega_{w}.
	\end{equation*} 
	For any $f\in H$, we have 
	\begin{align*}
		f&=UT_{P\Lambda P}^{\ast}f\\&=U\bigg(\int_{\Omega}\int_{\Omega_{w}}\langle \Lambda_{w}Pf,e_{w,v}\rangle e_{w,v}\delta_{w}\ d\mu(v) \ d\mu(w)\bigg)\\&=\int_{\Omega}\int_{\Omega_{w}}\langle f,P\Lambda_{w}^{\ast}e_{w,v}\rangle U(e_{w,v}\delta_{w})\ d\mu(v)\ d\mu(w)\\&=\int_{\Omega}\int_{\Omega_{w}}\langle f,Pu_{w,v}\rangle Q\Gamma_{w}^{\ast}e_{w,v}\ d\mu(v)\ d\mu(w)\\&=\int_{\Omega}Q\Gamma_{w}^{\ast}\bigg(\int_{\Omega_{w}}\langle Pf,u_{w,v}\rangle e_{w,v}\ d\mu(v)\bigg)\ d\mu(w)\\&=\int_{\Omega}Q\Gamma_{w}^{\ast}\Lambda_{w}Pf\ d\mu(w),
	\end{align*}
	where $u_{w,v}=\Lambda_{w}^{\ast}e_{w,v}$. By the definition of  controlled continuous dual $g$-frame,  $\{\Gamma_{w}\}_{w}$ is a  $(P,Q)$-controlled continuous dual $g$-frame of $\{\Lambda_{w}\}_{w}$.
	
	On the other hand, suppose that a  $(Q,Q)$-controlled continuous $g$-frame $\{\Gamma_{w}\}_{w}$ is a  $(P,Q)$-controlled continuous dual $g$-frame of $\{\Lambda_{w}\}_{w}$. For any $f\in H$, we have 
	\begin{equation*}
		f=\int_{\Omega}P\Lambda_{w}^{\ast}\Gamma_{w}^{\ast}Qf\ d\mu(w)=\int_{\Omega}Q\Gamma_{w}^{\ast}\Lambda_{w}Pf\ d\mu(w),
	\end{equation*} 
	that is, $T_{Q\Lambda Q}T_{P\Lambda P}^{\ast}=I_{H}$. Let $U=T_{Q\Gamma Q}$. Then $U:\ell^{2}(\{H_{w}\}_{w})\rightarrow H$ is a bounded left-inverse of $T_{P\Lambda P}^{\ast}$. A calculation as above shows that 
	\begin{equation*}
		\int_{\Omega}\int_{\Omega_{w}}\langle f,Pu_{w,v}\rangle Q\Gamma_{w}^{\ast}e_{w,v}\ d\mu(v)\ d\mu(w)=f=\int_{\Omega}\int_{\Omega_{w}}\langle f,Pu_{w,v}\rangle U(e_{w,v}\delta_{w})\ d\mu(v)\ d\mu(w), \forall f\in H.
	\end{equation*}
	Combining this with the fact $\{e_{w,v}\}_{v\in \Omega_{w}}$ is a continuous orthonormal basis of $H_{w}$, we have 
	\begin{equation*}
		Q\Gamma_{w}^{\ast}e_{w,v}=U(e_{w,v}\delta_{w}),\quad w\in\Omega, v\in\Omega_{w}.
	\end{equation*}
\end{proof}
\begin{theorem}
	Let $P\in GL^{+}(H)$ and $\{\Lambda_{w}\}_{w}$ be a  $(P,P)$-controlled continuous $g$-frame for $H$ with respect to $\{H_{w}\}_{w}$ with the synthesis operator and frame operator $T_{P\Lambda P}$ and $S_{P\Lambda P}$, respectively. Then $\{\Gamma_{w}\}_{w}$ is a  $P$-controlled continuous dual $g$-frame of $\{\Lambda_{w}\}_{w}$ if and only if 
	\begin{equation*}
		\Gamma_{w}f=(Tf)_{w}+\Lambda_{w}S_{P\Lambda P}^{-1}Pf, w\in \Omega,v\in H,
	\end{equation*}
	where $T:H\rightarrow \ell^{2}(\{H_{w}\}_{w})$ is a bounded linear operator satisfying $T_{P\Lambda P}T=0$.
\end{theorem}
\begin{proof}
	If  $T:H\rightarrow \ell^{2}(\{H_{w}\}_{w})$ is a bounded linear operator satisfying $T_{P\Lambda P}T=0$, then $\{\Gamma_{w}\}_{w}$ is a $g$-Bessel family for $H$ with respect to $\{H_{w}\}_{w}$. In fact, any $f\in H$ we have 
	\begin{align*}
		\int_{\Omega}\|\Gamma_{w}f\|^{2}\ d\mu(w)&=\int_{\Omega}\|(Tf)_{w}+\Lambda_{w}S_{P\Lambda P}^{-1}Pf\|^{2}\ d\mu(w)\\&\leq 2\bigg(\int_{\Omega}\|\Lambda_{w}S_{P\Lambda P}^{-1}Pf\|^{2}\ d\mu(w)+\|Tf\|^{2}\bigg)\\&\leq 2(B\|S_{P\Lambda P}^{-1}P\|^{2}+\|T\|^{2})\|f\|^{2},
	\end{align*}
	where $B$ is the upper bound of $\{\Lambda_{w}\}_{w\in\Omega}$. Furthermore,
	\begin{align*}
		\int_{\Omega}P\Lambda_{w}^{\ast}\Gamma_{w}f\ d\mu(w)&=\int_{\Omega}P\Lambda_{w}^{\ast}((Tf)_{w}+\Lambda_{w}S_{P\Lambda P}^{-1}Pf)\ d\mu(w)\\&=T_{P\Lambda P}Tf+\int_{\Omega}P\Lambda_{w}^{\ast}\Lambda_{w}S_{P\Lambda P}^{-1}Pf\ d\mu(w)=f.
	\end{align*}
	Thus $\{\Gamma_{w}\}_{w\in\Omega}$ is a  $P$-controlled continuous dual $g$-frame of $\{\Lambda_{w}\}_{w\in \Omega}$.
	
	Now we prove the converse. Assume that $\{\Gamma_{w}\}_{w\in\Omega}$ is a  $P$-controlled continuous dual $g$-frame of $\{\Lambda_{w}\}_{w\in \Omega}$. Define the operator $T$ as follows:
	\begin{equation*}
		T:H\rightarrow \ell^{2}(\{H_{w}\}_{w\in\Omega}),\quad f\rightarrow (Tf)_{w}, \quad f\in H
	\end{equation*}
	satisfying
	\begin{equation*}
		\Gamma_{w}f=(Tf)_{w}+\Lambda_{w}S_{P\Lambda P}^{-1}Pf, w\in\Omega.
	\end{equation*}
	For any $f\in H$, we have 
	\begin{align*}
		\|Tf\|^{2}&=\int_{\Omega}\|\Gamma_{w}f-\Lambda_{w}S^{-1}_{P\Lambda P}Pf \|^{2}\ d\mu(w)\\&\leq \int_{\Omega}\|\Gamma_{w}f\|^{2}\ d\mu(w)+\int_{\Omega}\|\Lambda_{w}S_{P\Lambda P}^{-1}Pf\|^{2}\ d\mu(w)\\&\quad+2\bigg(\int_{\Omega}\|\Gamma_{w}f\|^{2}\ d\mu(w)\bigg)^{\frac{1}{2}}\bigg(\int_{\Omega}\|\Lambda_{w}S_{P\Lambda P}^{-1}Pf\|^{2}\ d\mu(w)\bigg)^{\frac{1}{2}}\\&\leq (B_{1}+A^{-1}+2\sqrt{B_{1}A^{-1}})\|f\|^{2},
	\end{align*}
	where $B_{1}$ is the frame upper bound of $\{\Gamma_{w}\}_{w\in\Omega}$ and $A$ is the frame lower bound of $\{\Lambda_{w}\}_{w\in\Omega}$. Thus $T$ is a linear bounded operator. Moreover, for any $f,g\in H$, we have 
	\begin{align*}
		\langle T_{P\Lambda P}Tf,g\rangle&=\int_{\Omega}\langle P\Lambda_{w}^{\ast}Tf,g\rangle \ d\mu(w) \\&=\int_{\Omega}\langle P\Lambda_{w}^{\ast}(\Gamma_{w}f-\Lambda_{w}S_{P\Lambda P}^{-1}Pf),g\rangle \ d\mu(w) \\&=\int_{\Omega}\langle P\Lambda_{w}^{\ast}\Gamma_{w}f,g\rangle \ d\mu(w)-\int_{\Omega}\langle P\Lambda_{w}^{\ast}\Lambda_{w}S_{P\Lambda P}^{-1}Pf,g\rangle \ d\mu(w) \\&=\langle f,g\rangle -\langle f,g\rangle=0. 
	\end{align*}
	That is $T_{P\Lambda P}T=0$. The proof is completed.
\end{proof}
\medskip

\section*{Declarations}

\medskip

\noindent \textbf{Availablity of data and materials}\newline
\noindent Not applicable.

\medskip

\noindent \textbf{Human and animal rights}\newline
\noindent We would like to mention that this article does not contain any studies
with animals and does not involve any studies over human being.

\medskip

\noindent \textbf{Conflict of interest}\newline
\noindent The authors declare that they have no competing interests.

\medskip
%
\bibliographystyle{amsplain}

\begin{thebibliography}{99}
	
	\bibitem{AFC} M. R. Abdollahpour, M. H. Faroughi, Continuous G-frames in Hilbert spaces, Southeast Asian Bull. Math. 32 (2008), 1-19.
	
	\bibitem{YM} Y. Alizadeh,  A. Mohammad Reza,  Controlled  continuous  G-frames  and  their  multipliers in Hilbert spaces,  Sahand Communications in Mathematical Analysis (SCMA) Vol. 15 No. 1 (2019), 37-48. 
	
	\bibitem{COH} S. T. Ali, J. P. Antoine, J.P. Gazeau, Coherent states, Wavelets, and Their Generalizations, Springer-Verlag, Berlin, 2000.
	
	\bibitem{GAZ} S. T. Ali, J. P. Antoine, J. P. Gazeau,  Continuous frames in Hilbert spaces, Ann.  Phys. 
	222 (1993), 1-37.
	
	\bibitem{AHDR} A. Askari-Hemmat, M. A. Dehghan, M. Radjabalipour, \emph{Generalized frames and their redundancy}, Proc. Am. Math. Soc. {\bf 129} (2001), 1143-1147.
	
	\bibitem{PDA} P. Balazs, D. Bayer, A. Rahimi, Multipliers for continuous frames in Hilbert spaces, J. Phys. A Math. Theor. 45, 2240023(20p) (2012).
	
	\bibitem{CO}  O. Christensen, Frames and Bases: An Introductory Course. Birkhäuser, Boston. 2008.
	
	\bibitem{DGM} I. Daubechies, A. Grossmann,  Y. Meyer, Painless nonorthogonal expansions, J. Math.
	Phys. {\bf 27}  (1986), 1271-1283.	
	
	\bibitem{Duf} R. J. Duffin, A. C. Schaeffer, A class of nonharmonic fourier series, Trans. Am. Math. Soc. {\bf 72} (1952),
	341-366.
	
	\bibitem{GHN}  J. P. Gabardo,  D. Han,  \emph{Frames associated with measurable space}, Adv. Comput. Math. {\bf 18}
	(2003), no. 3, 127-147.	
	
	
	\bibitem{Gabor} D. Gabor.
	Theory of Communication, 
	J. IEE, 93(26):429-457, November 1946.
	
	
	\bibitem{KAI} G. Kaiser,  A Friendly Guide to Wavelets, Birkh\"{a}user, Boston, 1994.
	
	\bibitem{RM} M. Rahmani, Characterization of continuous g-frames via operators, Proceedings of the Institute of Mathematics and Mechanics,
	National Academy of Sciences of Azerbaijan
	Volume 46, Number 1, 2020, Pages 79-93.
	
	
\end{thebibliography}

\end{document}